\documentclass[11pt]{article}

\usepackage[T1]{fontenc}
\usepackage{lmodern}
\usepackage[margin=1in]{geometry}
\usepackage{amsmath,amssymb,amsthm,mathtools}
\usepackage{enumitem}
\usepackage{booktabs,array}
\usepackage{xcolor}
\usepackage[protrusion=true,expansion=false]{microtype}
\usepackage{fvextra}
\usepackage{xurl}
\usepackage{authblk}

\usepackage{tikz}
\usetikzlibrary{arrows.meta,fit,positioning}

\usepackage[colorlinks=true,linkcolor=blue,citecolor=blue,urlcolor=blue]{hyperref}
\hypersetup{
  pdftitle={A Lean Formalization of the Hamilton--Perelman Proof of the Three-Dimensional Poincare Conjecture},
  pdfauthor={Ziyang Qin, Yuan Liao, Ayush Khaitan, and Bennett Chow}
}

\theoremstyle{plain}
\newtheorem{theorem}{Theorem}

\theoremstyle{definition}

\theoremstyle{remark}

\newcommand{\Ric}{\operatorname{Ric}}

\newcommand{\releasetag}{v0.1.3}
\newcommand{\releasecommit}{7a48598d35109aa99d1cc678e2724c213cdf4ff3}
\newcommand{\releaseblobbase}{https://github.com/qinz1yang/differential-geometry/blob/\releasecommit/DifferentialGeometry/}
\newcommand{\releasetreebase}{https://github.com/qinz1yang/differential-geometry/tree/\releasecommit/DifferentialGeometry/}
\newcommand{\declareleanfile}[3]{%
  \expandafter\def\csname leanfile@#1\endcsname{#2}%
  \expandafter\def\csname leanline@#1\endcsname{#3}}
\newcommand{\leanref}[1]{%
  \ifcsname leanfile@#1\endcsname
    \href{\releaseblobbase\csname leanfile@#1\endcsname\#L\csname leanline@#1\endcsname}{\nolinkurl{#1}}%
  \else
    \nolinkurl{#1}%
  \fi}
\newcommand{\leanrefas}[2]{%
  \ifcsname leanfile@#1\endcsname
    \href{\releaseblobbase\csname leanfile@#1\endcsname\#L\csname leanline@#1\endcsname}{#2}%
  \else
    #2%
  \fi}
\newcommand{\leanfile}[1]{\nolinkurl{#1}}

\declareleanfile{exists_isManifold_three}{Topology/PiecewiseLinear/Moise352Producer.lean}{33}
\declareleanfile{plApproximation_three}{Topology/PiecewiseLinear/Moise352Producer.lean}{27}
\declareleanfile{plSmoothingCompact_three}{Topology/PiecewiseLinear/Moise352Producer.lean}{30}
\declareleanfile{smoothPoincareConjecture_holds}{Geometry/Flow/RicciFlow/Surgery/Skeleton/PoincareEndgame.lean}{103}
\declareleanfile{poincare_conjecture}{Topology/ThreeManifold/Poincare.lean}{12}
\declareleanfile{ricci_flow_short_time_existence}{Geometry/Flow/RicciFlow/ShortTime/Existence.lean}{34}
\declareleanfile{fixed_kappa_compactness}{Geometry/Flow/RicciFlow/Perelman/CanonicalNeighborhood/HighCurvatureModelBounds.lean}{566}
\declareleanfile{smooth_canonical_neighborhood}{Geometry/Flow/RicciFlow/Perelman/CanonicalNeighborhood/HighCurvatureModelBounds.lean}{539}
\declareleanfile{exists_poincare_controlled_extinction}{Geometry/Flow/RicciFlow/Surgery/Skeleton/PoincareEndgame.lean}{91}
\declareleanfile{isPoincareStandard_of_controlledExtinction}{Geometry/Flow/RicciFlow/Surgery/Topology/ExtinctionReconstruction.lean}{39}
\declareleanfile{fundamentalGroup_finiteConnectedSum_freeProduct}{Topology/VanKampen/FiniteConnectedSumFreeProduct.lean}{277}
\declareleanfile{positiveFreeContractibleClass}{Geometry/Flow/RicciFlow/Surgery/Topology/FreeLoopClass.lean}{1747}
\declareleanfile{canonicalWidth}{Geometry/Flow/RicciFlow/Extinction/Width/CanonicalClass.lean}{19}
\declareleanfile{historyWidth}{Geometry/Flow/RicciFlow/Extinction/Width/SurgeryWidthEvolution.lean}{385}
\declareleanfile{extinctionThreshold}{Geometry/Flow/RicciFlow/Extinction/Families/ScalarThreshold.lean}{11}
\makeatletter
\g@addto@macro{\UrlBreaks}{\do\_\do\.\do\/}
\makeatother
\title{A Lean Formalization of the Hamilton--Perelman Proof of the Three-Dimensional Poincar\'{e} Conjecture\footnote{Ziyang Qin and Yuan Liao contributed equally.}}
\author[1]{Ziyang Qin}
\author[2]{Yuan Liao}
\author[3]{Ayush Khaitan}
\author[2]{Bennett Chow}

\affil[1]{Department of Mathematics, Cornell University}
\affil[2]{Department of Mathematics, University of California San Diego}
\affil[3]{Princeton Language and Intelligence, DaIS, Princeton University}

\date{}

\begin{document}
\maketitle

We formalize the smooth three-dimensional Poincar\'{e}
conjecture, together with the Moise smoothing theorem,
yielding the topological three-dimensional Poincar\'{e}
conjecture.
The smooth proof follows the
Hamilton--Perelman route through Ricci flow with surgery and finite-time
extinction. Verification refers to the release in
Section~\ref{sec:verification}.

\section{Main theorems}

\begin{theorem}[Moise: existence of a smooth structure {\cite{Moise1977}}]
\label{thm:moise-smooth-existence}
Every closed topological three-manifold admits a \(C^\infty\) smooth
structure compatible with its given topology.
\end{theorem}

The Lean endpoint is
\leanrefas{exists_isManifold_three}{\mbox{\texttt{exists\_isManifold\_three}}}.
Throughout, \(S^3\subset\mathbb R^4\) carries its standard smooth structure.

\begin{theorem}[Smooth Poincar\'{e} conjecture; Hamilton--Perelman]
\label{thm:smooth-poincare}
Every simply-connected closed smooth three-manifold is diffeomorphic to
\(S^3\).
\end{theorem}

The Lean endpoint is \leanref{smoothPoincareConjecture_holds}.

\begin{theorem}[Topological Poincar\'{e} conjecture]
\label{thm:topological-poincare}
Every simply-connected closed topological three-manifold is homeomorphic
to \(S^3\).
\end{theorem}

The declaration \leanref{poincare_conjecture} applies Moise smoothing,
invokes the smooth theorem, and takes the underlying homeomorphism.
Its conclusion concerns the original topology. Its fully qualified type is:
\begin{Verbatim}[fontsize=\footnotesize,breaklines=true,commandchars=\\\{\}]
theorem DifferentialGeometry.Topology.poincare_conjecture
    (M : Type u) [TopologicalSpace M]
    [ChartedSpace (EuclideanSpace \ensuremath{\mathbb{R}} (Fin 3)) M]
    [T2Space M] [CompactSpace M] [SimplyConnectedSpace M] :
    Nonempty (M \ensuremath{\simeq\sb{\mathrm{t}}} Metric.sphere
      (0 : EuclideanSpace \ensuremath{\mathbb{R}} (Fin 4)) 1)
\end{Verbatim}
This proves the
\href{https://github.com/leanprover-community/mathlib4/blob/0df444a360eaa60ab8c11dca51a86af692955474/Mathlib/Geometry/Manifold/PoincareConjecture.lean\#L47}{three-dimensional topological Poincar\'{e} statement}
marked \texttt{proof\_wanted} in Mathlib~v4.33.1. Our smooth endpoint
likewise proves the
\href{https://github.com/leanprover-community/mathlib4/blob/0df444a360eaa60ab8c11dca51a86af692955474/Mathlib/Geometry/Manifold/PoincareConjecture.lean\#L52}{corresponding smooth statement}.
Both matches were checked in Lean by applying our endpoints directly,
with no additional hypotheses.

\section{Mathematical route}

For background on Ricci flow, see \cite{ChowEtAl2010,ChowEtAl2015}.
For accounts of Perelman's work, see Kleiner--Lott \cite{KleinerLott2008}
and Morgan--Tian \cite{MorganTian2007}, with the correction
\cite{MorganTian2015}.

Hamilton's earlier work includes convergence for three-manifolds with
positive Ricci curvature \cite{Hamilton1982} and for four-manifolds with
positive curvature operator \cite{Hamilton1986}, Ricci flow on surfaces
\cite{Hamilton1988}, the differential Harnack estimate
\cite{Hamilton1993Harnack}, the analysis of singularity formation
\cite{Hamilton1995Singularities}, and surgery for four-manifolds with
positive isotropic curvature~\cite{Hamilton1997}.

Starting with a smooth metric, the
\leanrefas{ricci_flow_short_time_existence}{short-time existence theorem}
uses the DeTurck method \cite{DeTurck1983,Hamilton1982} to construct a solution of
\[\partial_t g=-2\Ric.\]
Hamilton--Ivey pinching \cite{Hamilton1995Singularities,Ivey1993},
Perelman noncollapsing \cite{Perelman2002}, and pointed compactness
\cite{Hamilton1995Compactness}
control the singularity models. The main interfaces are
\leanrefas{fixed_kappa_compactness}{compactness of normalized ancient
\(\kappa\)-solutions} with fixed \(\kappa>0\) and
\leanrefas{smooth_canonical_neighborhood}{Perelman's canonical neighborhood
theorem}. Their surgery counterparts \cite{Perelman2003Surgery} control
neck cuts, cap attachments, and restarts, with a positive volume cost per
cut on a fixed time horizon.

For extinction, take a closed oriented simply-connected smooth
three-manifold \((M,g_0)\). We use the following least-area width. For
such a manifold \(N\), let \(\Lambda_0N\) be the space of contractible
continuous loops with the compact-open topology. Poincar\'{e} duality gives
\(H_2(N;\mathbb Z)\cong H^1(N;\mathbb Z)=0\), so Hurewicz gives
\(\pi_2(N)=0\). The Hurewicz theorem then identifies \(\pi_3(N)\) with
\(H_3(N;\mathbb Z)\). The oriented fundamental class, via this
identification and loop adjunction, determines a nontrivial free homotopy
class \(\alpha_N\) of maps \(S^2\to\Lambda_0N\), formalized as
\leanref{positiveFreeContractibleClass}. Since \(\pi_2(N)=0\), this class
is essential: no family of constant loops represents \(\alpha_N\). Define
\[
 A_g(\gamma)=\inf_{u|_{\partial D^2}=\gamma}\operatorname{Area}_g(u),
 \qquad
 W(N,g)=\inf_{[\Gamma]=\alpha_N}\sup_{z\in S^2}A_g(\Gamma(z)).
\]
Here \(u:D^2\to N\) ranges over Lipschitz filling disks, and \(\Gamma\)
ranges over continuous \(S^2\)-families of contractible \(C^1\) loops,
with the \(C^1\) topology induced by their values and first jets; the
homotopy class is taken in \(\Lambda_0N\).
This finite, nonnegative quantity is \leanref{canonicalWidth}.

Choose \(c>0\) from the initial metric so that
\(R_{g_0}\geq-3/(2c)\). The scalar-curvature lower barrier is preserved
through surgery, so \(R_{g(t)}\geq-3/(2(t+c))\) on every surviving
component. We identify the initial stage with \((M,g_0)\) by an
orientation-preserving isometry. For any component surviving to a finite
observation time \(T\), trace its ancestors through the surgery history
and let \(W(t)\) be their width at time \(t\), as in
\leanref{historyWidth}. The surviving components remain simply connected,
and every such ancestor chain starts at \(M\); hence
\(W(0)=W_0:=W(M,g_0)\) is independent of the chosen terminal component.
Let \(R_{\min}(t)\) denote the minimum scalar curvature on the ancestor
component defining \(W(t)\).

The least-area argument \cite{Hamilton1999,Perelman2003Extinction} gives,
at regular times \(0\leq t<T\),
\[
 D^+W(t)\leq-2\pi-\tfrac12 R_{\min}(t)W(t)
          \leq-2\pi+\frac{3}{4(t+c)}W(t),
\]
where \(D^+W(t):=\limsup_{h\downarrow0}(W(t+h)-W(t))/h\) is the upper
right Dini derivative. The width is continuous between surgeries and
right-continuous on \([0,T)\). At each surgery time \(s\), degree-one comparison maps preserve
the chosen loop-space class; their Lipschitz constants tend to one as
the incoming time approaches \(s\). Thus
\(W(s)\leq\liminf_{t\uparrow s}W(t)\).
The integrating factor \((t+c)^{-3/4}\), together with this jump
inequality, yields
\[
 0\leq\frac{W(t)}{(t+c)^{3/4}}
 \leq\frac{W_0}{c^{3/4}}
       -8\pi\bigl((t+c)^{1/4}-c^{1/4}\bigr).
\]
Consequently no component can survive beyond the common bound
\[
 T_*:=\left(c^{1/4}+\frac{W_0}{8\pi c^{3/4}}\right)^4-c,
\]
which is \leanref{extinctionThreshold}. Together with canonical-neighborhood
continuation and finiteness of surgery events on fixed time horizons,
obtained using the uniform positive volume cost per cut, this proves the
\leanrefas{exists_poincare_controlled_extinction}{controlled finite-time
extinction theorem}: every simply-connected closed oriented smooth
three-manifold, with any initial smooth metric, admits a finite surgery
history ending in the empty manifold at a positive finite time.

The history also records that every discarded closed connected component
is diffeomorphic to a connected sum of spherical space forms and
\(S^2\times S^1\) factors. Using this control,
\leanrefas{isPoincareStandard_of_controlledExtinction}{cut-and-cap
reconstruction} expresses the initial manifold as a connected sum of the
same types of factors. The
\leanrefas{fundamentalGroup_finiteConnectedSum_freeProduct}{van Kampen
theorem} identifies its fundamental group as the free product of the
factor groups. Simple connectivity excludes \(S^2\times S^1\) and
nontrivial spherical quotients; the remaining sphere factors yield
\(S^3\). Finally,
\leanrefas{plApproximation_three}{three-dimensional PL approximation} and
\leanrefas{plSmoothingCompact_three}{compact PL smoothing} supply the
compatible smooth atlas for the topological theorem. A fuller account of
these constructions will appear elsewhere.

\section{Formal artifact and verification}
\label{sec:verification}

The formalization uses Lean~4 \cite{Lean2021} and Mathlib
\cite{Mathlib2020}.
As in \cite{ChowLiaoQin2026}, we fix the source revision, here the
\href{https://github.com/qinz1yang/differential-geometry/releases/tag/\releasetag}{release \releasetag}:

\begin{center}
\small
\begin{tabular}{@{}ll@{}}
\toprule
Repository & \url{https://github.com/qinz1yang/differential-geometry}\\
Tag & \texttt{\releasetag}\\
Commit & \href{https://github.com/qinz1yang/differential-geometry/commit/\releasecommit}{\texttt{\releasecommit}}\\
Lean / Mathlib & \texttt{v4.33.1} / \texttt{v4.33.1}\\
\bottomrule
\end{tabular}
\end{center}

All links to our source tree use the full release commit SHA.
The toolchain and package dependencies are pinned by
\texttt{lean-toolchain} and \texttt{lake-manifest.json}.
Vendored sources and local adaptations are fixed by the release commit,
with upstream revisions and modification records preserved alongside the
sources.

The library incorporates and adapts selected developments from other Lean
projects: the DeGiorgi elliptic-regularity library of Scott Armstrong and
Julia Kempe \cite{ArmstrongKempe2026,DeGiorgiSoftware};
\'{A}lvaro Begu\'{e}'s Jordan--Sch\"{o}nflies development
\cite{SchoenfliesSoftware}; planar topology and triangulation constructions
from the SF LEAN meetup's Classification of Compact Surfaces project
\cite{ClassificationOfSurfacesSoftware}; and selected Sard, Morse-genericity,
and immersion-chart arguments from Tau Ceti \cite{TauCetiSoftware}.
These uses concern selected source developments, not the full theorem
collections of ClassificationOfSurfaces or Tau Ceti.
The Riemann-mapping adaptations come from Yury Kudryashov's specific
\texttt{urkud/mathlib4} development revision associated with Mathlib
PR~\#33505 \cite{RiemannMappingSoftware}, separately from the pinned
Mathlib package.

Our derivative-homotopy code adapts the local homotopy construction in
\texttt{Solution.lean}, lines 62156--62228, from \texttt{plby/HopfProblem}
\cite{HopfProblemSoftware}; this attribution concerns that construction,
not the entire Hopf formalization.
The preserved source maps and modification records distinguish these
origins and local extensions, including code moved into this project's
namespaces. We further acknowledge Jack McCarthy's contribution to
\href{\releaseblobbase Bundle/Equiv.lean}{bundle equivalences}, and Heather
Macbeth's and Jack McCarthy's contributions to
\href{\releaseblobbase Tensor/Alternating/Bundle/Defs.lean}{alternating-map bundles},
as credited in the preserved source headers. Both modules lie in the
import closure of the topological Poincar\'{e} endpoint.

\begin{samepage}
From the repository root at this release, build the endpoint and the
complete library with:

\begin{Verbatim}[fontsize=\footnotesize,breaklines=true]
lake build DifferentialGeometry.Topology.ThreeManifold.Poincare
lake build DifferentialGeometry
\end{Verbatim}
\end{samepage}

To print the transitive axiom closures, append the following commands to
\begin{center}
\small
\href{\releaseblobbase Topology/ThreeManifold/Poincare.lean}{\texttt{DifferentialGeometry/Topology/ThreeManifold/Poincare.lean}}
\end{center}
and rerun the first build command above:

\begin{Verbatim}[fontsize=\footnotesize,breaklines=true]
open DifferentialGeometry.Topology.PiecewiseLinear
open DifferentialGeometry.PDE.RicciFlow.Surgery.Topology

#print axioms DifferentialGeometry.Topology.poincare_conjecture
#print axioms exists_isManifold_three
#print axioms smoothPoincareConjecture_holds
#print axioms exists_poincare_controlled_extinction
\end{Verbatim}

These four queries return the same axiom set:
\[
 [\texttt{propext},\ \texttt{Classical.choice},\ \texttt{Quot.sound}].
\]
The queried proof closures contain neither \texttt{sorryAx} nor a
project-specific axiom.

These checks concern the named theorem closures. Their mathematical
interpretation requires review of the definitions and hypotheses.
Automated assistance was used for proofs and prose; the authors retain
responsibility for the mathematics and exposition.

\section{Acknowledgements}
A.K.~acknowledges funding from the DARPA expMath Program, as well as the AMS-Simons Travel Grant.

\end{document}